\documentclass[12pt]{article}
\usepackage{graphicx,amsmath,amssymb,amsfonts,times,hyperref}

\begin{document}

\title{Human Operator Modeling and Lie-Derivative Based Control}
\author{Tijana T. Ivancevic and Bojan Jovanovic}
\date{}
\maketitle\vspace{3cm}

\abstract{The motivation behind mathematically modeling the
\textit{human operator} is to help explain the response
characteristics of the complex dynamical system including the
human manual controller. In this paper, we present two approaches to
human operator modeling: classical linear control approach and
modern nonlinear control approach. The latter one is formalized using both fixed and adaptive Lie-Derivative based controllers.\bigbreak\bigbreak

}\newpage

\tableofcontents

\section{Introduction}

Despite the increasing trend toward automation, robotics and
artificial intelligence (AI) in many environments, the \emph{human
operator} will probably continue for some time to be integrally
involved in the control and regulation of various machines (e.g.,
missile--launchers, ground vehicles, watercrafts, submarines,
spacecrafts, helicopters, jet fighters, etc.). A typical manual
control task is
the task in which control of these machines is accomplished by \textit{%
manipulation of the hands or fingers} \cite{Wickens}. As
human--computer interfaces evolve, interaction techniques
increasingly involve a much more continuous form of interaction
with the user, over both human--to--computer (input) and
computer--to--human (output) channels. Such interaction could
involve gestures, speech and animation in addition to more
`conventional' interaction via mouse, joystick and keyboard. This
poses a problem for the design of interactive systems as it
becomes increasingly necessary to consider interactions occurring
over an interval, in continuous time.

The so--called \textit{manual control theory} developed out of the efforts
of feedback control engineers during and after the World War II, who
required models of human performance for continuous military tasks, such as
tracking with anti--aircraft guns \cite{Wiener}. This seems to be an area
worth exploring, firstly since it is generally concerned with systems which
are controlled in continuous time by the user, although discrete time
analogues of the various models exist. Secondly, it is an approach which
models both system and user and hence is compatible with research efforts on
`synthetic' models, in which aspects of both system and user are specified
within the same framework. Thirdly, it is an approach where continuous
mathematics is used to describe functions of time. Finally, it is a theory
which has been validated with respect to experimental data and applied
extensively within the military domains such as avionics.

The premise of manual control theory is that for certain tasks, the
performance of the human operator can be well approximated by a describing
function, much as an inanimate controller would be. Hence, in the literature
frequency domain representations of behavior in continuous time are applied.
Two of the main classes of system modelled by the theory are \textit{%
compensatory} and \textit{pursuit} systems. A system where only the error
signal is available to the human operator is a compensatory system. A system
where both the target and current output are available is called a pursuit
system. In many pursuit systems the user can also see a portion of the input
in advance; such tasks are called \textit{preview tasks} \cite{doh}.

A simple and widely used model is the `crossover model' \cite{2},
which has two main parameters, a \textit{gain} $K$ and a
\textit{time delay} $\tau$, given by the transfer function in the
Laplace transform $s$ domain
$$H=K{\mathrm{e}^{-\tau s}\over s}.
$$
Even with this simple model we can investigate some quite
interesting phenomena. For example consider a compensatory system
with a certain delay, if we have a low gain, then the system will
move only slowly towards the target, and hence will seem sluggish.
An expanded version of the crossover model is given by the
transfer function \cite{Wickens}
$$ H=K\frac{(T_Ls+1)\,\mathrm{e}^{-(\tau s+\alpha/s)}}{(T_Is+1)(T_Ns+1)},
$$
where $T_L$ and $T_I$ are the lead and lag constants (which
describe the \textit{equalization} of the human operator), while
the first--order lag $(T_NS+1)$ approximates the neuromuscular lag
of the hand and arm. The expanded term $\alpha/s$ in the time
delay accounts for the `phase drop', i.e., increased lags observed
at very low frequency \cite{GaneshSprBig}.

Alternatively if the gain $K$ is very high, then the system is
very likely to overshoot the target, requiring an adjustment in
the opposite direction, which may in turn overshoot, and so on.
This is known as `oscillatory behavior'. Many more detailed models
have also been developed; there are `anthropomorphic models',
which have a cognitive or physiological basis. For example the
`structural model' attempts to reflect the structure of the human,
with central nervous system, neuromuscular and vestibular
components \cite{doh}. Alternatively there is the `optimal control
modeling' approach, where algorithmic models which very closely
match empirical data are used, but which do not have any direct
relationship or explanation in terms of human neural and cognitive
architecture \cite{3}. In this model, an operator is assumed to
perceive a vector of displayed quantities and must exercise
control to minimize a \textit{cost functional} given by
\cite{Wickens}
$$J=E\{\lim_{T\to\infty} \frac{1}{T} \int^T_0 [q_iy_i^2(t) + \sum_i(r_iu^2(t)
+ g_i\dot{u}^2(t))]dt \},
$$
which means that the operator will attempt to minimize  the
expected value $E$ of some weighted combination of squared display
error $y$, squared control displacement $u$ and squared control
velocity $\dot u$. The relative values of the weighting constants
$q_i, r_i, g_i$ will depend upon the relative importance of
control precision, control effort and fuel expenditure.

In the case of manual control of a vehicle, this modeling yields
the `closed--loop' or `operator--vehicle' dynamics. A quantitative
explanation of this closed--loop behavior is necessary to
summarize operator behavioral data, to understand operator control
actions, and to predict the operator--vehicle dynamic
characteristics. For these reasons, control engineering
methodologies are applied to modeling human operators. These
`control theoretic' models primarily attempt to represent the
operator's control behavior, not the physiological and
psychological structure of the operator \cite{0,00}. These models
`gain in acceptability' if they can identify features of these
structures, `although they cannot be rejected' for failing to do
so \cite{1}.

One broad division of human operator models is whether they
simulated a continuous or discontinuous operator control strategy.
Significant success has been achieved in modeling human operators
performing compensatory and pursuit tracking tasks by employing
continuous, quasi--linear operator models. Examples of these
include the crossover optimal control models mentioned above.

Discontinuous input behavior is often observed during manual control of
large amplitude and acquisition tasks \cite{2,4,5,6}. These discontinuous
human operator responses are usually associated with \textit{precognitive
human control behavior} \cite{2,7}. Discontinuous control strategies have
been previously described by `bang--bang' or relay control techniques. In
\cite{8}, the authors highlighted operator's preference for this type of
relay control strategy in a study that compared controlling high--order
system plants with a proportional verses a relay control stick. By allowing
the operator to generate a sharper step input, the relay control stick
improved the operators' performance by up to 50 percent. These authors
hypothesized that when a human controls a high--order plant, the operator
must consider the error of the system to be dependent upon the integral of
the control input. Pulse and step inputs would reduce the integration
requirements on the operator and should make the system error response more
predictable to the operator.

Although operators may employ a \emph{bang--bang control} strategy, they often
impose an internal limit on the magnitude of control inputs. This internal
limit is typically less than the full control authority available \cite{2}.
Some authors \cite{9} hypothesized that this behavior is due to the
operator's recognition of their own reaction time delay. The operator must
tradeoff the cost of a switching time error with the cost of limiting the
velocity of the output to a value less than the maximum.

A significant amount of research during the 1960's and 1970's examined
discontinuous input behavior by human operators and developed models to
emulate it \cite{7,10,11,12,13,14,15,16,17}. Good summaries of these efforts
can be found in \cite{18}, \cite{4}, \cite{2} and \cite{0,00}. All of these
efforts employed some type of \textit{relay} element to model the
discontinuous input behavior. During the 1980's and 1990's, pilot models
were developed that included switching or discrete changes in pilot behavior
\cite{19,20,21,22,5,6}.

Recently, the so-called `variable structure control' techniques were applied
to model human operator behavior during acquisition tasks \cite{0,00}. The
result was a coupled, multi--input model replicating the discontinuous
control strategy. In this formulation, a switching surface was the
mathematical representation of the human operator's control strategy. The
performance of the variable strategy model was evaluated by considering the
longitudinal control of an aircraft during the visual landing task.

In this paper, we present two approaches to
human operator modeling: classical linear control approach and
modern nonlinear Lie-Derivative based control approach.

\section{Classical Control Theory versus Nonlinear Dynamics and Control}

In this section we review classical feedback control theory
(see e.g., \cite{GaneshSprSml,GaneshSprBig,GCompl}) and contrast
it with nonlinear and stochastic dynamics (see e.g., \cite{TacaNODY,StrAttr,Complexity}).

\subsection{Basics of Kalman's Linear State--Space Theory}

Linear multiple input--multiple output (MIMO) control systems can
always be put into Kalman canonical state--space form of order
$n$, with $m$ inputs and $k$ outputs. In the case of
\textit{continual time} systems we have
state and output equation of the form%
\begin{eqnarray}
d\mathbf{x}/dt
&=&\mathbf{A}(t)\,\mathbf{x}(t)+\mathbf{B}(t)\,\mathbf{u}(t),
\label{st_eq_con} \\
\mathbf{y}(t)
&=&\mathbf{C}(t)\,\mathbf{x}(t)+\mathbf{D}(t)\,\mathbf{u}(t),
\nonumber
\end{eqnarray}%
while in case of \textit{discrete time} systems we have state and
output
equation of the form%
\begin{eqnarray}
\mathbf{x}(n+1)
&=&\mathbf{A}(n)\,\mathbf{x}(n)+\mathbf{B}(n)\,\mathbf{u}(n),
\label{st_eq_dis} \\
\mathbf{y}(n)
&=&\mathbf{C}(n)\,\mathbf{x}(n)+\mathbf{D}(n)\,\mathbf{u}(n).
\nonumber
\end{eqnarray}%
Both in (\ref{st_eq_con}) and in (\ref{st_eq_dis}) the variables
have the following meaning:

$\mathbf{x}(t)\in \mathbb{X}$ is an $n-$vector of \textit{state
variables} belonging to the \textit{state space}
$\mathbb{X}\subset \mathbb{R}^{n}$;

$\mathbf{u}(t)\in \mathbb{U}$ is an $m-$vector of \textit{inputs}
belonging to the \textit{input space} $\mathbb{U}\subset
\mathbb{R}^{m}$;

$\mathbf{y}(t)\in \mathbb{Y}$ is a $k-$vector of \textit{outputs}
belonging to the \textit{output space} $\mathbb{Y}\subset
\mathbb{R}^{k}$;

$\mathbf{A}(t):$ $\mathbb{X}\rightarrow \mathbb{X}$ is an $n\times
n$ matrix of \textit{state dynamics};

$\mathbf{B}(t):$ $\mathbb{U}\rightarrow \mathbb{X}$ is an $n\times
m$ matrix of \textit{input map};

$\mathbf{C}(t):$ $\mathbb{X}\rightarrow \mathbb{Y}$ is an $k\times
n$ matrix of \textit{output map};

$\mathbf{D}(t):$ $\mathbb{U}\rightarrow \mathbb{Y}$ is an $k\times
m$ matrix of \textit{input--output transform}.

\bigbreak

Input $\mathbf{u}(t)\in \mathbb{U}$ can be empirically determined
by trial and error; it is properly defined by optimization process
called \textit{Kalman regulator}, or more generally (in the
presence of noise), by \textit{Kalman filter} (even better,
\textit{extended Kalman filter} to deal with stochastic
nonlinearities).

\subsection{Linear Stationary Systems and Operators}

\label{linstate}

The most common special case of the general Kalman model
(\ref{st_eq_con}), with constant state, input and output matrices
(and relaxed boldface vector--matrix notation), is the so--called
\textit{stationary linear model}
\begin{equation}
\dot{x}=Ax+Bu,\qquad y=Cx.  \label{slm}
\end{equation}%
The stationary linear system (\ref{slm}) defines a variety of
operators, in particular those related to the following problems:
\begin{enumerate}
  \item regulators,
  \item end point controls,
  \item servomechanisms, and
  \item repetitive modes (see \cite{Brockett}).
\end{enumerate}

\subsubsection{Regulator Problem and the Steady State Operator}

Consider a variable, or set of variables, associated with a
dynamical system. They are to be maintained at some desired values
in the face of changing circumstances. There exist a second set of
parameters that can be adjusted so as to achieve the desired
regulation. The effecting variables are usually called
\textit{inputs} and the affected variables called
\textit{outputs}. Specific examples include the regulation of the
thrust of a jet engine by controlling the flow of fuel, as well as
the regulation of the oxygen content of the blood using the
respiratory rate.

Now, there is the steady state operator of particular relevance
for the regulator problem. It is
\begin{equation*}
y_{\infty }=-CA^{-1}Bu_{\infty },
\end{equation*}
which describes the map from constant values of $u$ to the
equilibrium value of $y$. It is defined whenever $A$ is invertible
but the steady state value
will only be achieved by a real system if, in addition, the eigenvalues of $%
A $ have negative real parts. Only when the rank of $CA^{-1}B$
equals the dimension of $y$ can we steer $y$ to an arbitrary
steady state value and hold it there with a constant $u$. A
nonlinear version of this problem plays a central role in robotics
where it is called the \textit{inverse kinematics problem} (see,
e.g., \cite{Murray}).

\subsubsection{End Point Control Problem and the Adjustment Operator}

Here we have inputs, outputs and trajectories. In this case the
shape of the trajectory is not of great concern but rather it is
the end point that is of primary importance. Standard examples
include rendezvous problems such as one has in space exploration.

Now, the operator of relevance for the end point control problem,
is the operator
\begin{equation*}
x(T)=\int_{0}^{T}\exp [A(T-\sigma )]\,Bu(\sigma )\,d\sigma .
\end{equation*}%
If we consider this to define a map from the $m$D $L_{2}$ space $%
L_{2}^{m}[0,T]$ (where $u$ takes on its values) into
$\mathbb{R}^{m}$ then, if it is an onto map, it has a
Moore--Penrose (least squares) inverse
\begin{equation*}
u(\sigma )=B^{T}\exp [A^{T}(T-\sigma )]\,\left( W[0,T]\right)
^{-1}\left( x(T)-\exp (AT)\,x(0)\right) ,
\end{equation*}%
with the symmetric positive definite matrix $W$, the
\textit{controllability Gramian}, being given by
\begin{equation*}
W[0,T]=\int_{0}^{T}\exp [A(T-\sigma )]\,BB^{T}\exp [A^{T}(T-\sigma
)]\,d\sigma .
\end{equation*}

\subsubsection{Servomechanism Problem and the Corresponding Operator}

Here we have inputs, outputs and trajectories, as above, and an
associated dynamical system. In this case, however, it is desired
to cause the outputs to follow a trajectory specified by the
input. For example, the control of an airplane so that it will
travel along the flight path specified by the flight controller.

Now, because we have assumed that $A$, $B$ and $C$ are constant
\begin{equation*}
y(t)=C\exp (At)\,x(0)+\int_{0}^{t}C\exp [A(T-\tau )]\,Bu(\tau
)\,d\tau ,
\end{equation*}
and, as usual, the Laplace transform can be used to convert
convolution to multiplication. This brings out the significance of
the Laplace transform pair

\begin{equation}
C\exp (At)B\Longleftrightarrow C(Is-A)^{-1}B  \label{lapl}
\end{equation}
as a means of characterizing the input--output map of a linear
model with constant coefficients.

\subsubsection{Repetitive Mode Problem and the Corresponding Operator}

Here again one has some variable, or set of variables, associated
with a dynamical system and some inputs which influence its
evolution. The task has elements which are repetitive and are to
be done efficiently. Examples from biology include the control of
respiratory processes, control of the pumping action of the heart,
control of successive trials in practicing a athletic event.

The relevant operator is similar to the servomechanism operator,
however the constraint that $u$ and $x$ are periodic means that
the relevant diagonalization is provided by Fourier series, rather
than the Laplace transform. Thus, in the Fourier domain, we are
interested in a set of complex matrices
\begin{equation*}
G(iw_{i})=C(iw_{i}-A)^{-1}B,\qquad w_{i}=0,w_{0},2w_{0},...
\end{equation*}

More general, but still deterministic, models of the
input--state--output relation are afforded by the
\textit{nonlinear affine model} (see, e.g., \cite{Isidori})
\begin{eqnarray*}
\dot{x}(t) &=&f(x(t))+g(x(t))\,u(t), \\
y(t) &=&h(x(t));
\end{eqnarray*}
and the still more general \textit{fully nonlinear model}
\begin{eqnarray*}
\dot{x}(t) &=&f(x(t),u(t)), \\
y(t) &=&h(x(t)).
\end{eqnarray*}

\subsubsection{Feedback Changes the Operator}

No idea is more central to automatic control than the idea of
feedback. When an input is altered on the basis of the difference
between the actual output of the system and the desired output,
the system is said to involve \textit{feedback}. Man made systems
are often constructed by starting with a basic element such as a
motor, a burner, a grinder, etc. and then adding sensors and the
hardware necessary to use the measurement generated by the sensors
to regulate the performance of the basic element. This is the
essence of \textit{feedback control}. Feedback is often contrasted
with open loop systems in which the inputs to the basic element is
determined without reference to any measurement of the
trajectories. When the word feedback is used to describe naturally
occurring systems, it is usually implicit that the behavior of the
system can best be explained by pretending that it was designed as
one sees man made systems being designed \cite{Brockett}.

In the context of linear systems, the effect of feedback is easily
described. If we start with the stationary linear system
(\ref{slm}) with $u$ being the controls and $y$ being the measured
quantities, then the effect of feedback is to replace $u$ by
$u-Ky$ with $K$ being a matrix of feedback gains. The closed--loop
equations are then
\begin{equation*}
\dot{x}=(A-BKC)\,x+Bu,\qquad y=Cx.
\end{equation*}
Expressed in terms of the Laplace transform pairs (\ref{lapl}),
feedback effects the transformation
\begin{equation*}
\left( C\exp (At)B;C(Is-A)^{-1}B\right) \longmapsto C\exp
(A-BKC)^{t}B;C(Is-A+BKC)^{-1}B.
\end{equation*}
Using such a transformation, it is possible to alter the dynamics
of a system in a significant way. The modifications one can effect
by feedback include influencing the location of the eigenvalues
and consequently the stability of the system. In fact, if $K$ is
$m$ by $p$ and if we wish to select a gain matrix $K$ so that
$A-BKC$ has eigenvalues $\lambda _{1},\lambda _{2},...,\lambda
_{n}$, it is necessary to insure that
\begin{equation*}
\det \left(
\begin{array}{cc}
C(I\lambda _{1}-A)^{-1}B & -I \\
I & K%
\end{array}
\right) =0,\qquad i=1,2,...,n.
\end{equation*}

Now, if $CB$ is invertible then we can use the relationship $C\dot{x}%
=CAx+CBu $ together with $y=Cx$ to write $\dot{y}=CAx+CBu$. This
lets us solve for $u$ and recast the system as
\begin{eqnarray*}
\dot{x} &=&(A-B(CB)^{-1}CA)\,x+B(CB)^{-1}\dot{y}, \\
u &=&(CB)^{-1}\dot{y}-(CB)^{-1}CA\,x.
\end{eqnarray*}
Here we have a set of equations in which the roles of $u$ and $y$
are reversed. They show how a choice of $y$ determines $x$ and how
$x$ determines $u$ \cite{Brockett}.

\subsection{Stability and Boundedness}

Let a time--varying dynamical system may be expressed as
\begin{equation}
\dot{x}(t)=f(t,x(t))  \label{c43}
\end{equation}%
where $x\in \mathbb{R}^{n}$ is an $n$D vector and
$f:\mathbb{R}^{+}\times
D\rightarrow \mathbb{R}^{n}$ with $D=\mathbb{R}^{n}$ or $D=B_{h}$ for some $%
h>0$, where $B_{h}=\{x\in \mathbb{R}^{n}:\left\vert x\right\vert
<h\}$ is a ball centered at the origin with a radius of $h$. If
$D=\mathbb{R}^{n}$ then we say that the dynamics of the system are
defined \textit{globally}, whereas if $D=B_{h}$ they are only
defined \textit{locally}. We do not consider systems whose
dynamics are defined over disjoint subspaces of R. It is assumed
that $f(t,x)$ is piecemeal continuous in $t$ and Lipschitz in $x$
for existence and uniqueness of state solutions. As an example,
the linear system $\dot{x}(t)=Ax(t)$ fits the form of (\ref{c43})
with $D=\mathbb{R}^{n} $ \cite{Spooner}.

Assume that for every $x_{0}$ the initial value problem
\begin{equation*}
\dot{x}(t)=f(t,x(t)),\qquad x(t_{0})=x_{0},
\end{equation*}
possesses a unique solution $x(t,t_{0},x_{0});$ it is called a
solution to (\ref{c43}) if $x(t,t_{0},x_{0})=x_{0}$ and $\frac{d}{dt}%
x(t,t_{0},x_{0})=f(t,x(t,t_{0},x_{0}))$ \cite{Spooner}.

A point $x_{e}\in \mathbb{R}^{n}$ is called an \textit{equilibrium
point} of (\ref{c43}) if $f(t,x_{e})=0$ for all $t\geq 0$. An
equilibrium point $x_{e}$ is called an \textit{isolated
equilibrium point} if there exists an $\rho >0$
such that the ball around $x_{e}$, $B_{\rho }(x_{e})=\{x\in \mathbb{R}%
^{n}:\left\vert x-x_{e}\right\vert <\rho \},$ contains no other
equilibrium points besides $x_{e}$ \cite{Spooner}.

The equilibrium $x_{e}$ = 0 of (\ref{c43}) is said to be
\textit{stable in the sense of Lyapunov} if for every $\epsilon
>0$ and any $t_{0}\geq 0$ there exists a $\delta (\epsilon
,t_{0})>0$ such that $\left\vert
x(t,t_{0},x_{0})\right\vert <\epsilon $ for all $t\geq t_{0}$ whenever $%
\left\vert x_{0}\right\vert <\delta (\epsilon ,t_{0})$ and $%
x(t,t_{0},x_{0})\in B_{h}(x_{e})$ for some $h>0$. That is, the
equilibrium is stable if when the system (\ref{c43}) starts close
to $x_{e}$, then it will stay close to it. Note that stability is
a property of an equilibrium, not a system. A system is stable if
all its equilibrium points are stable. Stability in the sense of
Lyapunov is a local property. Also, notice that
the definition of stability is for a single equilibrium $x_{e}\in \mathbb{R}%
^{n}$ but actually such an equilibrium is a trajectory of points
that
satisfy the differential equation in (\ref{c43}). That is, the equilibrium $%
x_{e}$\ is a solution to the differential equation (\ref{c43}),
$x(t,t_{0},x_{0})=x_{e}$ for $t\geq 0$. We call any set such that
when the initial condition of (\ref{c43}) starts in the set and
stays in the set for all $t\geq 0$, an \textit{invariant set}. As
an example, if $x_{e}=0$ is an equilibrium, then the set
containing only the point $x_{e}$ is an invariant set, for
(\ref{c43}) \cite{Spooner}.

If $\delta $ is independent of $t_{0}$, that is, if $\delta
=\delta (\epsilon )$, then the equilibrium $x_{e}$ is said to be
\textit{uniformly
stable}. If in (\ref{c43}) $f$ does not depend on time (i.e., $f(x)$), then $%
x_{e}$ being stable is equivalent to it being uniformly stable.
Uniform stability is also a local property.

The equilibrium $x_{e}=0$ of (\ref{c43}) is said to be
\textit{asymptotically stable} if it is stable and for every
$t_{0}\geq 0$ there exists $\eta (t_{0})>0$ such that
$\lim_{t\rightarrow \infty }\left| x(t,t_{0},x_{0})\right| =0$
whenever $\left| x_{0}\right| <\eta (t_{0})$. That is, it is
asymptotically stable if when it starts close to the equilibrium
it will converge to it. Asymptotic stability is also a local
property. It is a stronger stability property since it requires
that the solutions to the ordinary differential equation converge
to zero in addition to what is required for stability in the sense
of Lyapunov.

The equilibrium $x_{e}=0$ of (\ref{c43}) is said to be
\textit{uniformly asymptotically stable} if it is uniformly stable
and for every $\epsilon >0$ and and $t_{0}\geq 0$, there exist a
$\delta _{0}>0$ independent of $t_{0}$
and $\epsilon $, and a $T(\epsilon )>0$ independent of $t_{0}$, such that $%
\left| x(t,t_{0},x_{0})-x_{e}\right| \leq \epsilon $ for all
$t\geq
t_{0}+T(\epsilon )$ whenever$\ \left| x_{0}-x_{e}\right| <\delta (\epsilon )$%
. Again, if in (\ref{c43}) $f$ does not depend on time (i.e., $f(x)$), then $%
x_{e}$ being asymptotically stable is equivalent to it being
uniformly asymptotically stable. Uniform asymptotic stability is
also a local property.

The set $X_{d}\subset \mathbb{R}^{n}$ of all $x_{0}\in
\mathbb{R}^{n}$ such that $\left\vert x(t,t_{0},x_{0})\right\vert
\rightarrow 0$ as $t\rightarrow
\infty $ is called the \textit{domain of attraction} of the equilibrium $%
x_{e}=0$ of (\ref{c43}). The equilibrium $x_{e}=0$ is said to be
\textit{asymptotically stable} in the large if $X_{d}\subset
\mathbb{R}^{n}$. That is, an equilibrium is asymptotically stable
in the large if no matter where the system starts, its state
converges to the equilibrium asymptotically. This is a global
property as opposed to the earlier stability definitions that
characterized local properties. This means that for asymptotic
stability in the large, the local property of asymptotic stability
holds for $B_{h}(x_{e})$ with $h=\infty $ (i.e., on the whole
state--space).

The equilibrium $x_{e}=0$ is said to be \textit{exponentially
stable} if
there exists an $\alpha >0$ and for every $\epsilon >0$ there exists a $%
\delta (\epsilon )>0$ such that $\left| x(t,t_{0},x_{0})\right|
\leq \epsilon e^{-\alpha (t-t_{0})},$ whenever $\left|
x_{0}\right| <\delta (\epsilon )$ and $t\geq t_{0}\geq 0$. The
constant $\alpha $ is sometimes called the \textit{rate of
convergence}. Exponential stability is sometimes said to be a
`stronger' form of stability since in its presence we know that
system trajectories decrease exponentially to zero. It is a local
property; here is its global version. The equilibrium point
$x_{e}=0$ is \textit{exponentially stable in the large} if there
exists $\alpha >0$ and for any $\beta >0$ there exists $\epsilon
(\beta )>0$ such that $\left| x(t,t_{0},x_{0})\right| \leq
\epsilon (\beta )e^{-\alpha (t-t_{0})},$ whenever $\left|
x_{0}\right| <\beta $ and $t\geq t_{0}\geq 0$.

An equilibrium that is not stable is called \textit{unstable}.

Closely related to stability is the concept of
\textit{boundedness}, which is, however, a global property of a
system in the sense that it applies to trajectories (solutions) of
the system that can be defined over all of the state--space
\cite{Spooner}.

A solution $x(t,t_{0},x_{0})$ of (\ref{c43}) is \textit{bounded}
if there exists a $\beta >0$, that may depend on each solution,
such that $\left\vert x(t,t_{0},x_{0})\right\vert <\beta $ for all
$t\geq t_{0}\geq 0$. A system is said to possess \textit{Lagrange
stability} if for each $t_{0}\geq 0$ and $x_{0}\in
\mathbb{R}^{n}$, the solution $x(t,t_{0},x_{0})$ is bounded. If an
equilibrium is asymptotically stable in the large or exponentially
stable in the large then the system for which the equilibrium is
defined is also Lagrange stable (but not necessarily vice versa).
Also, if an equilibrium is stable, it does not imply that the
system for which the equilibrium is defined is Lagrange stable
since there may be a way to pick $x_{0}$ such that it is near an
unstable equilibrium and $x(t,t_{0},x_{0})\rightarrow \infty $ as
$t\rightarrow \infty $.

The solutions $x(t,t_{0},x_{0})$ are \textit{uniformly bounded} if for any $%
\alpha >0$ and $t_{0}\geq 0$, there exists a $\beta (\alpha )>0$
(independent of $t_{0}$) such that if $\left| x_{0}\right| <\alpha $, then $%
\left| x(t,t_{0},x_{0})\right| <\beta (\alpha )$ for all $t\geq t_{0}\geq 0$%
. If the solutions are uniformly bounded then they are bounded and
the system is Lagrange stable.

The solutions $x(t,t_{0},x_{0})$ are said to be \textit{uniformly
ultimately bounded} if there exists some $B>0$, and if
corresponding to any $\alpha >0$ and $t_{0}>0$ there exists a
$T(\alpha )>0$ (independent of $t_{0}$) such that $\left|
x_{0}\right| <\alpha $ implies that $\left|
x(t,t_{0},x_{0})\right| <B$ for all $t\geq t_{0}+T(\alpha )$.
Hence, a system is said to be uniformly ultimately bounded if
eventually all trajectories end up in a $B-$neighborhood of the
origin.

\subsection{Lyapunov's Stability Method}

A. M. Lyapunov invented two methods to analyze stability
\cite{Spooner}. In his \textit{indirect method} he showed that if
we linearize a system about an equilibrium point, certain
conclusions about local stability properties can be made (e.g., if
the eigenvalues of the linearized system are in the left half
plane then the equilibrium is stable but if one is in the right
half plane it is unstable).

In his \textit{direct method} the stability results for an
equilibrium $x_{e}=0$ of (\ref{c43}) depend on the existence of an
appropriate \textit{Lyapunov function} $V:D\rightarrow \mathbb{R}$
where $D=\mathbb{R}^{n}$ for global results (e.g., asymptotic
stability in the large) and $D=B_{h}$ for some $h>0$, for local
results (e.g., stability in the sense of Lyapunov or asymptotic
stability). If $V$ is continuously differentiable with respect to
its arguments then the derivative of $V$ with respect to $t$ along
the solutions of (\ref{c43}) is
\begin{equation*}
\dot{V}(t,x)=\frac{\partial V}{\partial t}+\frac{\partial V}{\partial x}%
f(t,x).
\end{equation*}%
As an example, suppose that (\ref{c43}) is autonomous, and let
$V(x)$ is a quadratic form $V(x)=x^{T}Px$ where $x\in
\mathbb{R}^{n}$ and $P=P^{T}$.
Then, $\dot{V}(x)=\frac{\partial V}{\partial x}f(t,x)=\dot{x}^{T}Px+x^{T}P%
\dot{x}=2x^{T}P\dot{x}$ \cite{Spooner}.

Lyapunov's direct method provides for the following ways to test
for stability. The first two are strictly for local properties
while the last two have local and global versions.

- \emph{Stable}: If $V(t,x)$ is continuously differentiable,
positive definite, and $\dot{V}(t,x)\leq 0$, then $x_{e}=0$ is
stable.

- \emph{Uniformly stable}: If $V(t,x)$ is continuously
differentiable,
positive definite, decrescent\footnote{%
A $C^0-$function $V(t,x):\mathbb{R}^{+}\times B_{h}\rightarrow \mathbb{R%
}(V(t,x):\mathbb{R}^{+}\times \mathbb{R}^{n}\rightarrow
\mathbb{R})$ is said
to be \textit{decrescent} if there exists a strictly increasing function $%
\gamma $ defined on $[0,r)$ for some $r>0$ (defined on $[0,\infty
)$) such that $V(t,x)\leq \gamma (\left| x\right| )$ for all
$t\geq 0$ and $x\in B_{h} $ for some $h>0$.}, and $V(t,x)\leq 0$,
then $x_{e}=0$ is uniformly stable.

- \emph{Uniformly asymptotically stable}: If $V(t,x)$ is
continuously
differentiable, positive definite, and decrescent, with negative definite $%
\dot{V}(t,x)$, then $x_{e}=0$ is uniformly asymptotically stable
(uniformly asymptotically stable in the large if all these
properties hold globally).

- \emph{Exponentially stable}: If there exists a continuously
differentiable $V(t,x)$ and $c,c_{1},c_{2},c_{3}>0$ such that
\begin{eqnarray}
c_{1}\left\vert x\right\vert ^{c} &\leq &V(t,x)\leq
c_{2}\left\vert
x\right\vert ^{c},  \label{c56} \\
\dot{V}(t,x) &\leq &-c_{31}\left\vert x\right\vert ^{c},
\label{c57}
\end{eqnarray}%
for all $x\in B_{h}$ and $t\geq 0$, then $x_{e}=0$ is
exponentially stable. If there exists a continuously
differentiable function $V(t,x)$ and Equations (\ref{c56}) and
(\ref{c57}) hold for some $c,c_{1},c_{2},c_{3}>0$ for all $x\in
\mathbb{R}^{n}$ and $t\geq 0$, then $x_{e}=0$ is exponentially
stable in the large \cite{Spooner}.

\subsection{Nonlinear and Impulse Dynamics of Complex Plants}

In this section we give two examples of nonlinear dynamical
systems that are beyond reach of the classical control theory.

\subsubsection{Hybrid Dynamical Systems of Variable Structure}

Consider a \textit{hybrid dynamical system of variable structure},
given by $n-$dimensional ODE (see \cite{Michel})
\begin{equation}
\dot{x}=f(t,x),  \label{varsys}
\end{equation}
where $x=x(t)\in \mathbb{R}^{n}$ and $f=f(t,x):\mathbb{R}^{+}\times \mathbb{R%
}^{n}\rightarrow \mathbb{R}^{n}$. Let the domain $G\subset \mathbb{R}%
^{+}\times \mathbb{R}^{n},$ on which the vector--field $f(t,x)$ is
defined,
be divided into two subdomains, $G^{+}$ and $G^{-}$, by means of a smooth $%
(n-1)-$manifold $M$. In $G^{+}\cup M$, let there be given a vector--field $%
f^{+}(t,x)$, and in $G^{-}\cup M,$ let there be given a vector--field $%
f^{-}(t,x)$. Assume that both $f^{+}=f^{+}(t,x)$ and
$f^{-}=f^{-}(t,x)$ are continuous in $t$ and smooth in $x$. For
the system (\ref{varsys}), let
\begin{equation*}
f=\left\{
\begin{array}{c}
f^{+}\text{ \ \ when \ }x\in G^{+} \\
f^{-}\text{ \ \ when \ }x\in G^{-}%
\end{array}
\right. .
\end{equation*}
Under these conditions, a solution $x(t)$ of ODE (\ref{varsys}) is
well--defined while passing through $G$ until the manifold $M$ is
reached.

Upon reaching the manifold $M$, in physical systems with inertia,
the transition
\begin{equation*}
\text{from \ \ }\dot{x}=f^{-}(t,x)\text{ \ \ \ to \ \
}\dot{x}=f^{+}(t,x)
\end{equation*}
does not take place instantly on reaching $M$, but after some
delay. Due to this delay, the solution $x(t)$ oscillates about
$M$, $x(t)$ being displaced along $M$ with some mean velocity.

As the delay tends to zero, the limiting motion and velocity along
$M$ are determined by the \textit{linear homotopy ODE}
\begin{equation}
\dot{x}=f^{0}(t,x)\equiv (1-\alpha )\,f^{-}(t,x)+\alpha
\,f^{+}(t,x), \label{homsys}
\end{equation}
where $x\in M$ and $\alpha \in \lbrack 0,1]$ is such that the
\textit{linear
homotopy segment} $f^{0}(t,x)$ is tangential to $M$ at the point $x$, i.e., $%
f^{0}(t,x)\in T_{x}M$, where $T_{x}M$ is the tangent space to the manifold $%
M $ at the point $x$.

The vector--field $f^{0}(t,x)$ of the system (\ref{homsys}) can be
constructed as follows: at the point $x\in M,$ $f^{-}(t,x)$ and
$f^{+}(t,x)$ are given and their ends are joined by the linear
homotopy segment. The point of intersection between this segment
and $T_{x}M$ is the end of the required vector--field
$f^{0}(t,x)$. The vector function $x(t)$ which
satisfies (\ref{varsys}) in $G^{-}$ and $G^{+}$, and (\ref{homsys}) when $%
x\in M,$ can be considered as a \textit{solution} of
(\ref{varsys}) in a \textit{general sense}.

However, there are cases in which the solution $x(t)$ cannot
consist of a finite or even countable number of arcs, each of
which passes through $G^{-}$ or $G^{+}$ satisfying (\ref{varsys}),
or moves along the manifold $M$ and satisfies the homotopic ODE
(\ref{homsys}). To cover such cases, assume that the vector--field
$f=f(t,x)$ in ODE (\ref{varsys}) is a Lebesgue--measurable
function in a domain $G\subset \mathbb{R}^{+}\times
\mathbb{R}^{n}$, and that for any closed bounded domain $D\subset
G$ there exists a summable function $K(t)$ such that almost
everywhere in $D$ we have $|f(t,x)|\leq K(t) $. Then the
absolutely continuous vector function $x(t)$ is called the
\textit{generalized solution} of the ODE (\ref{varsys}) \textit{in
the sense
of Filippov} (see \cite{Michel}) if for almost all $t$, the vector $\dot{x}=%
\dot{x}(t)$ belongs to the least convex closed set containing all
the limiting values of the vector field $f(t,x^{\ast })$, where
$x^{\ast }$
tends towards $x$ in an arbitrary manner, and the values of the function $%
f(t,x^{\ast })$ on a set of measure zero in $\mathbb{R}^{n}$ are
ignored.

Such \textit{hybrid systems} of variable structure occur in the
study of nonlinear electric networks (endowed with electronic
switches, relays, diodes, rectifiers, etc.), in models of both
natural and artificial neural networks, as well as in feedback
control systems (usually with continuous--time plants and digital
controllers/filters).

\subsubsection{Impulse Dynamics of Kicks and Spikes}

\paragraph{The Spike Function.}

Recall that the \textit{Dirac's }$\delta -$\textit{function} (also
called the \textit{\
impulse function} in the systems and signals theory) represents a \textit{%
limit} of the \textit{Gaussian bell--shaped curve}
\begin{equation}
g(t,\alpha )=\frac{1}{\sqrt{\pi \alpha }}\mathrm{e}^{-t^{2}/\alpha
}\text{ \ \ (with parameter }\alpha \rightarrow 0\text{)}
\label{Gauss}
\end{equation}
() where the factor $1/\sqrt{\pi \alpha }$ serves for the \textit{%
normalization} of (\ref{Gauss}),
\begin{equation}
\int_{-\infty }^{+\infty }\frac{dt}{\sqrt{\pi \alpha }}\mathrm{e}%
^{-t^{2}/\alpha }=1,  \label{GaussNormal}
\end{equation}
i.e., the \textit{area} under the \textit{pulse} is equal to
\textit{unity}. In (\ref{Gauss}), the smaller $\alpha $\ the
higher the peak. In other words,
\begin{equation}
\delta (t)=\lim_{\alpha \rightarrow 0}\frac{1}{\sqrt{\pi \alpha }}\mathrm{e}%
^{-t^{2}/\alpha },  \label{delataGauss}
\end{equation}
which is a pulse so short that outside of $t=0$ it vanishes,
whereas at $t=0$ it still remains normalized according to
(\ref{GaussNormal}). Therefore, we get the usual definition of the
$\delta -$function:
\begin{eqnarray}
\delta (t) &=&0\text{ \ \ for \ \ t}\neq 0,  \notag \\
\int_{-\epsilon }^{+\epsilon }\delta (t)\,dt &=&1,  \label{delta}
\end{eqnarray}
where $\epsilon $ may be arbitrarily small. Instead of centering
the $\delta -$pulse around $t=0$, we can center it around any
other time $t_{0}$ so that (\ref{delta}) is transformed into
\begin{eqnarray}
\delta (t-t_{0}) &=&0\text{ \ \ for \ \ t}\neq t_{0},  \notag \\
\int_{t_{0}-\epsilon }^{t_{0}+\epsilon }\delta (t-t_{0})\,dt &=&1.
\label{shiftDelta}
\end{eqnarray}

Another well--known fact is that the integral of the $\delta
-$function is the \textit{Heaviside's step function}
\begin{equation}
H(T)=\int_{-\infty }^{T}\delta (t)\,dt=\left\{
\begin{array}{c}
\begin{array}{c}
0\text{ \ \ for \ \ }T<0 \\
1\text{ \ \ for \ \ }T>0
\end{array}
\\
(\frac{1}{2}\text{ \ \ for \ \ }T=0)
\end{array}
\right. .  \label{Heaviside}
\end{equation}
Now we can perform several generalizations of the relation (\ref{Heaviside}%
). First, we have
\begin{equation*}
\int_{-\infty }^{T}\delta (ct-t_{0})\,dt=\left\{
\begin{array}{c}
\begin{array}{c}
0\text{ \ \ for \ \ }T<t_{0}/c \\
1/c\text{ \ \ for \ \ }T>t_{0}/c \\
\frac{1}{2c}\text{ \ \ for \ \ }T=t_{0}/c
\end{array}
\end{array}
\right. .
\end{equation*}
More generally, we can introduce the so--called \textit{phase function} $%
\phi (t),$ (e.g., $\phi (t)=ct-t_{0}$) which is continuous at
$t=t_{0}$ but its time derivative $\dot{\phi}(t)\equiv \frac{d\phi
(t)}{dt}$ is discontinuous at $t=t_{0}$ (yet positive,
$\dot{\phi}(t)>0$), and such that
\begin{equation*}
\int_{-\infty }^{T}\delta (\phi (t))\,dt=\left\{
\begin{array}{c}
\begin{array}{c}
0\text{ \ \ for \ \ }T<t_{0} \\
1/\dot{\phi}(t_{0})\text{ \ \ for \ \ }T>t_{0} \\
\frac{1}{\dot{\phi}(t_{0})}\text{ \ \ for \ \ }T=t_{0}
\end{array}
\end{array}
\right. .
\end{equation*}
Finally, we come the the \textit{spike function }$\delta (\phi (t))\dot{\phi}%
(t)$, which like $\delta -$function represents a \textit{spike} at $t=t_{0},$%
\ such that the normalization criterion (\ref{shiftDelta}) is
still valid,
\begin{equation*}
\int_{t_{0}-\epsilon }^{t_{0}+\epsilon }\delta (\phi
(t))\dot{\phi}(t)\,dt=1.
\end{equation*}

\paragraph{Deterministic Delayed Kicks.}

Following Haken \cite{Haken}, we consider the mechanical example
of a soccer ball that is kicked by a soccer player and rolls over
grass, whereby its motion will be slowed down. In our opinion,
this is a perfect model for all `shooting--like' actions of the
human operator.

We start with the Newton's (second) law of motion,
$m\dot{v}=force$, and in order to get rid of superfluous
constants, we put temporarily $m=1$. The $force$ on the r.h.s.
consists of the damping force $-\gamma v(t) $ of the grass (where
$\gamma $ is the damping constant) and the sharp force
$F(t)=s\delta (t-\sigma )$ of the individual kick
occurring at time $t=\sigma$ (where $s$ is the strength of the kick, and $%
\delta$ is the Dirac's `delta' function). In this way, the
single--kick equation of the ball motion becomes
\begin{equation}
\dot{v}=-\gamma v(t)+s\delta (t-\sigma ),  \label{eqm}
\end{equation}%
with the general solution
\begin{equation*}
v(t)=sG(t-\sigma ),
\end{equation*}%
where $G(t-\sigma )$ is the Green's function\footnote{%
This is the Green's function of the first order system
(\ref{eqm}). Similarly, the Green's function
\begin{equation*}
G(t-\sigma )=\left\{
\begin{array}{c}
0\text{ \ \ \ for \ }t<\sigma \\
(t-\sigma )\mathrm{e}^{-\gamma (t-\sigma )}\text{ \ \ \ for \ }t\geq \sigma%
\end{array}%
\right.
\end{equation*}
corresponds to the second order system
\begin{equation*}
\left(\frac{d}{dt}+\gamma\right)^2 G(t-\sigma)=\delta(t-\sigma).
\end{equation*}%
}
\begin{equation*}
G(t-\sigma)=\left\{
\begin{array}{c}
0\text{ \ \ \ for \ }t<\sigma \\
\mathrm{e}^{-\gamma (t-\sigma )}\text{ \ \ \ for \ }t\geq \sigma%
\end{array}%
\right. .
\end{equation*}

Now, we can generalize the above to $N$ kicks with individual
strengths $s_j$, occurring at a sequence of times $\{\sigma
_{j}\}$, so that the total kicking force becomes
\begin{equation*}
F(t)=\sum_{j=1}^{N}s_{j}\delta (t-\sigma _{j}).
\end{equation*}
In this way, we get the multi--kick equation of the ball motion
\begin{equation*}
\dot{v}=-\gamma v(t)+\sum_{j=1}^{N}s_{j}\delta (t-\sigma _{j}),
\end{equation*}
with the general solution
\begin{equation}
v(t)=\sum_{j=1}^{N}s_{j}G(t-\sigma _{j}).  \label{dis}
\end{equation}

As a final generalization, we would imagine that the kicks are
continuously exerted on the ball, so that kicking force becomes
\begin{equation*}
F(t)=\int_{t_{0}}^{T}s(\sigma )\delta (t-\sigma )d\sigma \equiv
\int_{t_{0}}^{T}d\sigma F(\sigma )\delta (t-\sigma ),
\end{equation*}%
so that the continuous multi--kick equation of the ball motion becomes%
\begin{equation*}
\dot{v}=-\gamma v(t)+\int_{t_{0}}^{T}s(\sigma )\delta (t-\sigma
)d\sigma \equiv -\gamma v(t)+\int_{t_{0}}^{T}d\sigma F(\sigma
)\delta (t-\sigma ),
\end{equation*}%
with the general solution
\begin{equation}
v(t)=\int_{t_{0}}^{T}d\sigma F(\sigma )G(t-\sigma
)=\int_{t_{0}}^{T}d\sigma F(\sigma )\mathrm{e}^{-\gamma (t-\sigma
)}.  \label{con}
\end{equation}

\paragraph{Random Kicks and Langevin Equations.}

We now denote the times at which kicks occur by $t_{j}$ and
indicate their direction in a one--dimensional game by $(\pm
1)_{j}$, where the choice of the plus or minus sign is random
(e.g., throwing a coin). Thus the kicking force can be written in
the form
\begin{equation}
F(t)=s\sum_{j=1}^{N}\delta (t-t_{j})(\pm 1)_{j},  \label{force}
\end{equation}
where for simplicity we assume that all kicks have the same
strength $s$. When we observe many games, then we may perform an
average $<...>$\ over all these different \textit{performances},
\begin{equation}
<F(t)>=s<\sum_{j=1}^{N}\delta (t-t_{j})(\pm 1)_{j}>.  \label{av}
\end{equation}
Since the direction of the kicks is assumed to be independent of
the time at which the kicks happen, we may split (\ref{av}) into
the product
\begin{equation*}
<F(t)>=s<\sum_{j=1}^{N}\delta (t-t_{j})><(\pm 1)_{j}>.
\end{equation*}
As the kicks are assumed to happen with equal frequency in both
directions, we get the cancellation
\begin{equation*}
<(\pm 1)_{j}>=0,
\end{equation*}
which implies that the average kicking force also vanishes,
\begin{equation*}
<F(t)>=0.
\end{equation*}
In order to characterize the strength of the force (\ref{force}),
we
consider a quadratic expression in $F$, e.g., by calculating the \textit{%
correlation function} for two times $t,t^{\prime },$%
\begin{equation*}
<F(t)F(t^{\prime })>=s^{2}<\sum_{j}\delta (t-t_{j})(\pm
1)_{j}\sum_{k}\delta (t^{\prime }-t_{k})(\pm 1)_{k}>.
\end{equation*}
As the ones for $j\neq k$ will cancel each other and for $j=k$
will become $1 $, the correlation function becomes a single sum
\begin{equation}
<F(t)F(t^{\prime })>=s^{2}<\sum_{j}\delta (t-t_{j})\delta
(t^{\prime }-t_{k})>,  \label{cor}
\end{equation}
which is usually evaluated by assuming the \textit{Poisson
process} for the times of the kicks.

Now, proper description of random motion is given by
\textit{Langevin rate equation}, which describes the
\textit{Brownian motion}: when a particle is immersed in a fluid,
the velocity of this particle is slowed down by a force
proportional to its velocity and the particle undergoes a zig--zag
motion (the particle is steadily pushed by much smaller particles
of the liquid in a random way). In physical terminology, we deal
with the behavior of a system (particle) which is coupled to a
\textit{heat bath} or reservoir (namely the liquid). The heat bath
has two effects:

\begin{enumerate}
\item It decelerates the mean motion of the particle; and

\item It causes statistical fluctuation.
\end{enumerate}

The standard Langevin equation has the form
\begin{equation}
\dot{v}=-\gamma v(t)+F(t),  \label{Lang}
\end{equation}
where $F(t)$ is a \textit{fluctuating force} with the following
properties:

\begin{enumerate}
\item Its statistical average (\ref{av}) vanishes; and

\item Its correlation function (\ref{cor}) is given by
\begin{equation}
<F(t)F(t^{\prime })>=Q\delta (t-t_{0}),  \label{corLang}
\end{equation}%
where $t_{0}=T/N$ denotes the mean free time between kicks, and
 $Q=s^{2}/t_{0}$ is the \textit{random fluctuation}.
\end{enumerate}

The general solution of the Langevin equation (\ref{Lang}) is
given by (\ref{con}).

The average velocity vanishes, $<v(t)>=0$, as both directions are
possible and cancel each other. Using the integral solution
(\ref{con}) we get
\begin{equation*}
<v(t)v(t^{\prime })>=<\int_{t_{0}}^{t}d\sigma
\int_{t_{0}}^{t^{\prime }}d\sigma ^{\prime }F(\sigma )F(\sigma
^{\prime })\mathrm{e}^{-\gamma (t-\sigma )}\mathrm{e}^{-\gamma
(t^{\prime }-\sigma ^{\prime })}>,
\end{equation*}
which, in the steady--state, reduces to
\begin{equation*}
<v(t)v(t^{\prime })>=\frac{Q}{2\gamma }\mathrm{e}^{-\gamma
(t-\sigma )},
\end{equation*}
and for equal times
\begin{equation*}
<v(t)^{2}>=\frac{Q}{2\gamma }.
\end{equation*}

If we now repeat all the steps performed so far with $m\neq 1$,
the final result reads
\begin{equation}
<v(t)^{2}>=\frac{Q}{2\gamma m}.  \label{vt2}
\end{equation}

Now, according to thermodynamics, the \textit{mean kinetic energy}
of a particle is given by
\begin{equation}
\frac{m}{2}<v(t)^{2}>=\frac{1}{2}k_{B}T,  \label{vt2T}
\end{equation}%
where $T$ is the (absolute) temperature, and $k_{B}$ is the
Boltzman's constant. Comparing (\ref{vt2}) and (\ref{vt2T}), we
obtain the important Einstein's result
\begin{equation*}
Q=2\gamma k_{B}T,
\end{equation*}%
which says that whenever there is damping, i.e., $\gamma \neq 0,$
then there are random fluctuations (or noise) $Q$. In other words,
fluctuations or noise are inevitable in any physical system. For
example, in a resistor (with the resistance $R$) the electric
field $E$ fluctuates with a
correlation function (similar to (\ref{corLang}))%
\begin{equation*}
<E(t)E(t^{\prime })>=2Rk_{B}T\delta (t-t_{0}).
\end{equation*}%
This is the simplest example of the so--called
\textit{dissipation--fluctuation theorem}.

\section{Nonlinear Control Modeling of the\\ Human Operator}

In this section we present the basics of modern nonlinear control,
as a powerful tool for controlling nonlinear dynamical systems.

\subsection{Graphical Techniques for Nonlinear Systems}

\begin{figure}[tbh]
\centerline{\includegraphics[width=9cm]{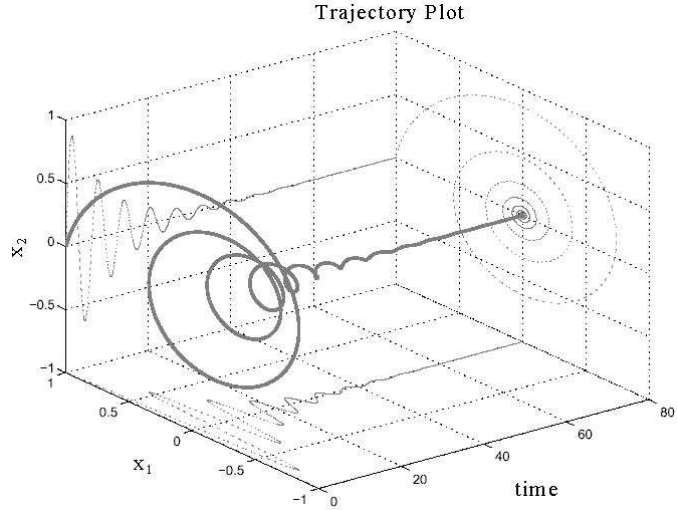}}
\caption{Tracer plot of the damped spring.} \label{TracerPlot}
\end{figure}

Graphical techniques preceded modern geometrical techniques in
nonlinear control theory. They started with simple plotting tools,
like the so--called `tracer plot'. It is a useful visualization
tool for analysis of second order dynamical systems, which just
adds time dimension to the standard 2D
phase portrait. For example, consider the damped spring governed by%
\begin{equation*}
\ddot{x}=-k\dot{x}-x,\ x(0)=1.
\end{equation*}%
Its tracer plot is given in Figure \ref{TracerPlot}. Note the
stable asymptote reached as $t\rightarrow \infty $.

The most important graphical technique is the so--called
describing function analysis.

\subsubsection{Describing Function Analysis}

Describing function analysis extends classical linear control
technique, frequency response analysis, for nonlinear systems
\cite{Wilson}. It is an \textit{approximate} graphical method
mainly used to predict \textit{limit cycles} in nonlinear ODEs.

\begin{figure}[tbh]
\centerline{\includegraphics[width=11cm]{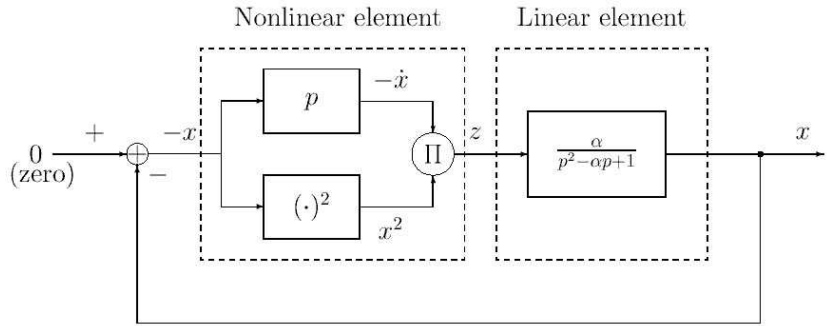}}
\caption{Feedback interpretation of the Van Der Pol oscillator
(after \protect\cite{Wilson}). Here $p$ is a (linear)
differentiator, and $\Pi $ a (nonlinear) multiplicator.}
\label{VDPblock}
\end{figure}

\begin{figure}[tbh]
\centerline{\includegraphics[width=10cm]{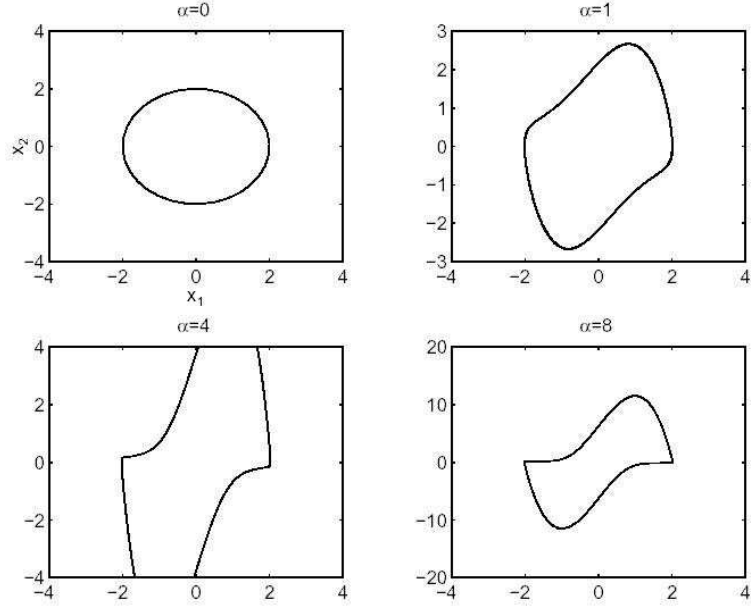}}
\caption{Limit cycle of the Van der Pol oscillator (after \protect\cite%
{Wilson}). Approximation is reasonable for small $\protect\alpha$,
but error, (amplitude should be equal to 2), grows as
$\protect\alpha$ is increased.} \label{VDPlimitC}
\end{figure}

For example, if we want to predict the existence of limit cycles
in the classical \textit{Van der Pol's oscillator} given by
\begin{equation}
\ddot{x}+\alpha (x^{2}-1)\,\dot{x}+x=0,  \label{vdp1}
\end{equation}%
we need to rewrite (\ref{vdp1}) as a \textit{linear} unstable
low--pass block and a \textit{nonlinear} block (see Figure
\ref{VDPblock}). In this
way, using the nonlinear block substitution, $w:=-\dot{x}x^{2}$, we get%
\begin{eqnarray*}
\ddot{x}-\alpha \dot{x}+x &=&\alpha w,\qquad\text{or} \\
x(p^{2}-\alpha p+1) &=&\alpha w,
\end{eqnarray*}%
or just considering the transfer function from $w$ to $x$,

\begin{equation*}
\frac{x}{w}=\frac{\alpha }{p^{2}-\alpha p+1}.
\end{equation*}%
Now, if we assume that the Van der Pol oscillator does have a
limit cycle with a frequency of
\begin{equation*}
x(t)=A\sin (wt),
\end{equation*}%
so $\dot{x}=Aw\cos (wt)$, therefore the \textit{output} of the
nonlinear block is
\begin{eqnarray*}
z &=&-x^{2}\dot{x}=-A^{2}\sin ^{2}(wt)\,Aw\cos (wt) \\
&=&-\frac{A^{3}w}{2}\left( 1-\cos (2wt)\right) \cos (wt) \\
&=&-\frac{A^{3}w}{4}\left( \cos (wt)-\cos (3wt)\right) .
\end{eqnarray*}%
Note how $z$\ contains a third harmonic term, but this is
attenuated by the low--pass nature of the linear block, and so
does not effect the signal in
the feedback. So we can approximate $z$ by%
\begin{equation*}
z\approx \frac{A^{3}}{4}w\cos
(wt)=\frac{A^{2}}{4}\frac{d}{dt}\left( -A\sin (wt)\right) .
\end{equation*}%
Therefore, the output of the nonlinear block can be approximated
by the quasi--linear transfer function which depends on the signal
amplitude, $A$, as well as frequency. The frequency response
function of the quasi--linear element is obtained by substituting
$p\equiv s=iw$,
\begin{equation*}
N(A,w)=\frac{A^{2}}{4}(iw).
\end{equation*}%
Since the system is assumed to contain a sinusoidal oscillation,
\begin{eqnarray*}
x &=&A\sin (wt)=G(iw)\,z \\
&=&G(iw)\,N(A,w)\,(-x),
\end{eqnarray*}%
where $G(iw)$ is the transfer function of the linear block. This
implies that,
\begin{equation*}
\frac{x}{-x}=-1=G(iw)\,N(A,w),
\end{equation*}%
so%
\begin{equation*}
1+\frac{A^{2}(iw)}{4}\frac{\alpha }{(iw)^{2}-\alpha (iw)+1}=0,
\end{equation*}%
which solving gives,
\begin{equation*}
A=2,\qquad \omega =1,
\end{equation*}%
which is independent of $\alpha $. Note that in terms of the
Laplace variable $p\equiv s$, the closed loop characteristic
equation of the system
is%
\begin{equation*}
1+\frac{A^{2}(iw)}{4}\frac{\alpha }{p^{2}-\alpha p+1}=0,
\end{equation*}%
whose eigenvalues are

\begin{equation*}
\lambda _{1,2}=-\frac{1}{8}\alpha (A^{2}-4)\pm \sqrt{\frac{\alpha
^{2}(A^{2}-4)^{2}}{64}-1.}
\end{equation*}%
Corresponding to $A=2$ gives eigenvalues of $\lambda _{1,2}=\pm i$
indicating an existence of a limit cycle of amplitude $2$ and
frequency $1$ (see Figure \ref{VDPlimitC}). If $A>2$ eigenvalues
are negative real, so stable, and the same holds for $A<2.$ The
approximation of the nonlinear block with $(A2/4)(iw)$ is called
the \textit{describing function}. This technique is useful because
most limit cycles are approximately sinusoidal and most linear
elements are low--pass in nature. So most of the higher harmonics,
if they existed, are attenuated and lost.

\subsection{Feedback Linearization}

The idea of feedback linearization is to algebraically transform
the nonlinear system dynamics into a fully or partly linear one so
that the linear control techniques can be applied. Note that this
is not the same as a conventional linearization using Jacobians.
In this subsection we will present the modern, geometrical,
Lie--derivative based techniques for exact feedback linearization
of nonlinear control systems.

\subsubsection{The Lie Derivative and Lie Bracket in Control Theory}
\label{liecon}

Recall that given a scalar function $h(x)$ and a
vector--field $f(x)$, we define a new scalar function, $\mathcal{L}%
_{f}h:=\nabla hf$, which is the \emph{Lie derivative} of $h$ w.r.t. $f$,
i.e., the directional derivative of $h$ along the direction of the
vector $f$ (see \cite{GaneshSprBig,GaneshADG}). Repeated Lie derivatives can be defined recursively:
\begin{eqnarray*}
\mathcal{L}_{f}^{0}h &=&h, \\
\mathcal{L}_{f}^{i}h &=&\mathcal{L}_{f}\left(
\mathcal{L}_{f}^{i-1}h\right) =\nabla \left(
\mathcal{L}_{f}^{i-1}h\right) f,\qquad \text{for \ }i=1,2,...
\end{eqnarray*}
Or given another vector--field, $g$, then
$\mathcal{L}_{g}\mathcal{L}_{f}h(x)$ is defined as
\begin{equation*}
\mathcal{L}_{g}\mathcal{L}_{f}h=\nabla \left(
\mathcal{L}_{f}h\right) g.
\end{equation*}

For example, if we have a control system
\begin{eqnarray*}
\dot{x} &=&f(x), \\
y &=&h(x),
\end{eqnarray*}
with the state $x=x(t)$ and the the output $y$, then the
derivatives of the output are:
\begin{eqnarray*}
\dot{y} &=&\frac{\partial h}{\partial
x}\dot{x}=\mathcal{L}_{f}h,\text{ \ \ \ \
and} \\
\ddot{y} &=&\frac{\partial L_{f}h}{\partial
x}\dot{x}=\mathcal{L}_{f}^{2}h.
\end{eqnarray*}

\begin{figure}[tbh]
\centerline{\includegraphics[height=4cm]{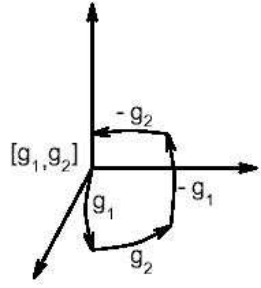}} \caption{`Lie
bracket motion' is possible by appropriately modulating the
control inputs (from \protect\cite{Goodwine}).} \label{liebrmot}
\end{figure}

Also, recall that the curvature of two vector--fields,
$g_{1},g_{2},$ gives a non--zero Lie bracket, $[g_{1},g_{2}]$ (see
Figure \ref{liebrmot}). Lie bracket motions can generate new
directions in which the system can move.

In general, the Lie bracket of two vector--fields, $f(x)$ and
$g(x),$ is defined by
\begin{equation*}
\left[ f,g\right] :=ad_{f}g:=\nabla gf-\nabla fg:=\frac{\partial
g}{\partial x}f-\frac{\partial f}{\partial x}g,
\end{equation*}%
where $\nabla f:=\partial f/\partial x$ is the Jacobian matrix. We
can define Lie brackets recursively,
\begin{eqnarray*}
ad_{f}^{0}g &=&g, \\
ad_{f}^{i}g &=&[f,ad_{f}^{i-1}g],\qquad \text{for \
}i=1,2,...\text{ }
\end{eqnarray*}%
Lie brackets have the properties of bilinearity,
skew--commutativity and Jacobi identity.

For example, if
\begin{equation*}
f=\left(
\begin{array}{c}
\cos x_{2} \\
x_{1}%
\end{array}%
\right) ,\qquad g=\left(
\begin{array}{c}
x_{1} \\
1%
\end{array}%
\right) ,
\end{equation*}%
then we have
\begin{eqnarray*}
\left[ f,g\right] &=&\left(
\begin{array}{cc}
1 & 0 \\
0 & 0%
\end{array}%
\right) \left(
\begin{array}{c}
\cos x_{2} \\
x_{1}%
\end{array}%
\right) -\left(
\begin{array}{cc}
0 & -\sin x_{2} \\
1 & 0%
\end{array}%
\right) \left(
\begin{array}{c}
x_{1} \\
1%
\end{array}%
\right) \\
&=&\left(
\begin{array}{c}
\cos x_{2}+\sin x_{2} \\
-x_{1}%
\end{array}%
\right).
\end{eqnarray*}

\subsubsection{Input/Output Linearization}

Given the single--input single--output (SISO) system
\begin{eqnarray}
\dot{x} &=&f(x)+g(x)\,u,  \label{siso} \\
y &=&h(x),  \notag
\end{eqnarray}
we want to formulate a linear differential equation relation between output $%
y$ and a new input $v$. We will investigate (see \cite{Isidori,Sastri,Wilson}%
):

\begin{itemize}
\item How to generate a linear input/output relation.

\item What are the internal dynamics and zero--dynamics associated
with the input/output linearization?

\item How to design stable controllers based on the I/O
linearization.
\end{itemize}

This linearization method will be exact in a finite domain, rather
than tangent as in the local linearization methods, which use
Taylor series approximation. Nonlinear controller design using the
technique is called exact feedback linearization.

\subsubsection{Algorithm for Exact Feedback Linearization}

\begin{figure}[tbh]
\centerline{\includegraphics[height=4cm]{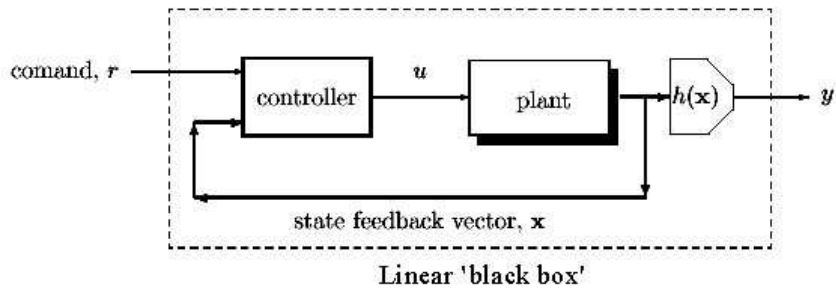}}
\caption{Feedback linearization.} \label{ExFebck}
\end{figure}

We want to find a nonlinear compensator such that the closed--loop
system is linear (see Figure \ref{ExFebck}). We will consider only
affine SISO systems of the type (\ref{siso}), i.e,
$\dot{x}=f(x)+g(x)\,u$, $y=h(x),$ and we will try to construct a
\textit{control law} of the form
\begin{equation}
u=p(x)+q(x)\,v,  \label{ctl}
\end{equation}
where $v$ is the setpoint, such that the \textit{closed--loop
system}
\begin{eqnarray*}
\dot{x} &=&f(x)+g(x)\,p(x)+g(x)\,q(x)\,v, \\
y &=&h(x),
\end{eqnarray*}
is linear from command $v$ to $y$.

The main idea behind the feedback linearization construction is to
find a nonlinear change of coordinates which transforms the
original system into one which is linear and controllable, in
particular, a chain of integrators. The difficulty is finding the
output function $h(x)$ which makes this construction possible.

We want to design an exact nonlinear feedback controller. Given
the nonlinear affine system, $\dot{x}=f(x)+g(x)$, $y=h(x)$, we
want to find the controller functions $p(x)$ and $q(x)$. The
unknown functions inside our controller (\ref{ctl}) are given by:
\begin{eqnarray}
p(x) &=&\frac{-\left( \mathcal{L}_{f}^{r}h(x)+\beta _{1}\mathcal{L}%
_{f}^{r-1}h(x)+...+\beta _{r-1}\mathcal{L}_{f}h(x)+\beta _{r}h(x)\right) }{%
\mathcal{L}_{g}\mathcal{L}_{f}^{r-1}h(x)},  \notag   \\
q(x) &=&\frac{1}{\mathcal{L}_{g}\mathcal{L}_{f}^{r-1}h(x)},
\label{pqx}
\end{eqnarray}
which are comprised of Lie derivatives, $\mathcal{L}_{f}h(x)$.
Here, the
\textit{relative order}, $r$, is the smallest integer $r$ such that $%
\mathcal{L}_{g}\mathcal{L}_{f}^{r-1}h(x)\neq 0$. For linear
systems $r$ is the difference between the number of poles and
zeros.

To obtain the \textit{desired response}, we choose the $r$
parameters in the $\beta $ polynomial to describe how the output
will respond to the setpoint, $v$ (pole--placement).
\begin{equation*}
\frac{d^{r}y}{dt^{r}}+\beta
_{1}\frac{d^{r-1}y}{dt^{r-1}}+...+\beta _{r-1} \frac{dy}{dt}+\beta
_{r}y=v.
\end{equation*}

Here is the proposed algorithm \cite{Isidori,Sastri,Wilson}):

\begin{enumerate}
\item Given nonlinear SISO process, $\dot{x}=f(x,u)$, and output equation $%
y=h(x)$, then:

\item Calculate the relative order, $r$.

\item Choose an $r$th order desired linear response using
pole--placement technique (i.e., select $\beta $). For this could
be used a simple $r$th order low--pass filter such as a
Butterworth filter.

\item Construct the exact linearized nonlinear controller
(\ref{pqx}), using Lie derivatives and perhaps a symbolic
manipulator (Mathematica or Maple).

\item Close the loop and obtain a linear input--output black--box
(see Figure \ref{ExFebck}).

\item Verify that the result is actually linear by comparing with
the desired response.
\end{enumerate}

\subsection{Controllability}

\subsubsection{Linear Controllability}

A system is \textit{controllable} if the set of all states it can
reach from
initial state $x_{0}=x(0)$ at the fixed time $t=T$ contains a ball $\mathcal{%
B}$ around $x_{0}$. Again, a system is \textit{small time locally
controllable} (STLC) iff the ball $\mathcal{B}$ for $\,t\leq T$
contains a
neighborhood of $x_{0}$.\footnote{%
The above definition of controllability tells us only whether or
not something can reach an open neighborhood of its starting
point, but does not
tell us how to do it. That is the point of the \textit{trajectory generation}%
.}

In the case of a linear system in the standard state--space form
\begin{equation}
\dot{x}=Ax+Bu,  \label{steq}
\end{equation}
where $A$ is the $n\times n$ \textit{state matrix} and $B$ is the
$m\times n$ \textit{input matrix}, all controllability definitions
coincide, i.e.,
\begin{eqnarray*}
0 &\rightarrow &x(T), \\
x(0) &\rightarrow &0, \\
x(0) &\rightarrow &x(T),
\end{eqnarray*}
where $T$ is either fixed or free.

\emph{Rank condition}\index{rank condition} states: System
(\ref{steq}) is controllable iff the matrix
\begin{equation*}
W_{n}=\left( B\,AB\,...\,A^{n-1}B\right) \qquad \text{has full
rank.}
\end{equation*}

In the case of nonlinear systems the corresponding result is
obtained using the formalism of Lie brackets, as Lie algebra is to
nonlinear systems as matrix algebra is to linear systems.

\subsubsection{Nonlinear Controllability}

Nonlinear MIMO--systems are generally described by differential
equations of the form (see \cite{Isidori,Schaft,Goodwine}):

\begin{equation}
\dot{x}=f(x)+g_{i}(x)\,u^{i},\qquad(i=1,...,n),  \label{NlSys}
\end{equation}
defined on a smooth $n-$manifold $M$, where $x\in M$ represents
the state of
the control system, $f(x)$ and $g_{i}(x)$ are vector--fields on $M$ and the $%
u^{i}$ are control inputs, which belong to a set of
\textit{admissible controls}, $u^{i}\in U$. The system
(\ref{NlSys}) is called \textit{driftless}, or \textit{kinematic},
or \textit{control linear} if $f(x)$ is identically zero;
otherwise, it is called a \textit{system with drift}, and the
vector--field $f(x)$ is called the \textit{drift term}. The flow
$\phi
_{t}^{g}(x_{0})$ represents the solution of the differential equation $\dot{x%
}=g(x)$ at time $t$ starting from $x_{0}$. Geometrical way to
understand the \textit{controllability} of the system
(\ref{NlSys}) is to understand the geometry of the vector--fields
$f(x)$ and $g_{i}(x)$.

\paragraph{Example: Car--Parking Using Lie Brackets.}
\label{carpark}

In this popular example, the driver has two different
transformations at his disposal. He can turn the steering wheel,
or he can drive the car forward or back. Here, we specify the
state of a car by four coordinates: the $(x,y)$ coordinates of the
center of the rear axle, the direction $\theta $ of the car, and
the angle $\phi $ between the front wheels and the direction of
the car. $L$ is the constant length of the car. Therefore, the
configuration manifold of the car is 4D, $M:=(x,y,\theta ,\phi )$.

Using (\ref{NlSys}), the driftless car kinematics can be defined
as:
\begin{equation}
\,\dot{x}=g_{1}(x)\,u_{1}+g_{2}(x)\,u_{2},  \label{Car2}
\end{equation}%
with two vector--fields $g_{1},g_{2}\in \mathcal{X}^{k}(M).$

The infinitesimal transformations will be the vector--fields
\begin{equation*}
g_{1}(x)\equiv \text{\textsc{drive}}=\cos \theta \frac{\partial }{\partial x}%
+\sin \theta \frac{\partial }{\partial y}+\frac{\tan \phi
}{L}\frac{\partial }{\partial \theta }\equiv \left(
\begin{array}{c}
\cos \theta \\
\sin \theta \\
\frac{1}{L}\tan \phi \\
0%
\end{array}
\right),
\end{equation*}
and
\begin{equation*}
g_{2}(x)\equiv \text{\textsc{steer}}=\frac{\partial }{\partial
\phi }\equiv \left(
\begin{array}{c}
0 \\
0 \\
0 \\
1%
\end{array}
\right) .
\end{equation*}

Now, \textsc{steer}\ and \textsc{drive}\ do not commute; otherwise
we could do all your steering at home before driving of on a trip.
Therefore, we have a Lie bracket
\begin{equation*}
\lbrack g_{2},g_{1}]\equiv \lbrack \text{\textsc{steer}},\text{\textsc{drive}%
}]=\frac{1}{L\cos ^{2}\phi }\frac{\partial }{\partial \theta }\equiv \text{%
\textsc{rotate}}.
\end{equation*}
The operation $[g_{2},g_{1}]\equiv $ \textsc{rotate }$\equiv \lbrack $%
\textsc{steer}$,$\textsc{drive}$]$ is the infinitesimal version of
the sequence of transformations: steer, drive, steer back, and
drive back, i.e.,
\begin{equation*}
\{\text{\textsc{steer}},\text{\textsc{drive}},\text{\textsc{steer}}^{-1},%
\text{\textsc{drive}}^{-1}\}.
\end{equation*}
Now, \textsc{rotate}\ can get us out of some parking spaces, but
not tight ones: we may not have enough room to \textsc{rotate}
out. The usual tight
parking space restricts the \textsc{drive}\ transformation, but not \textsc{%
steer}. A truly tight parking space restricts \textsc{steer}\ as
well by putting your front wheels against the curb.

Fortunately, there is still another commutator available:
\begin{eqnarray*}
\lbrack g_{1},[g_{2},g_{1}]] &\equiv &[\text{\textsc{drive}},[\text{\textsc{%
steer}},\text{\textsc{drive}}]]=[[g_{1},g_{2}],g_{1}]\equiv \\
\lbrack \text{\textsc{drive}},\text{\textsc{rotate}}]
&=&\frac{1}{L\cos
^{2}\phi }\left( \sin \theta \frac{\partial }{\partial x}-\cos \theta \frac{%
\partial }{\partial y}\right) \equiv \text{\textsc{slide}}.
\end{eqnarray*}
The operation $[[g_{1},g_{2}],g_{1}]\equiv $ \textsc{slide }$\equiv \lbrack $%
\textsc{drive}$,$\textsc{rotate}$]$ is a displacement at right
angles to the car, and can get us out of any parking place. We
just need to remember to steer, drive, steer back, drive some
more, steer, drive back, steer back, and drive back:
\begin{equation*}
\{\text{\textsc{steer}},\text{\textsc{drive}},\text{\textsc{steer}}^{-1},%
\text{\textsc{drive}},\text{\textsc{steer}},\text{\textsc{drive}}^{-1},\text{%
\textsc{steer}}^{-1},\text{\textsc{drive}}^{-1}\}.
\end{equation*}
We have to reverse steer in the middle of the parking place. This
is not intuitive, and no doubt is part of the problem with
parallel parking.

Thus from only two controls $u_{1}$ and $u_{2}$ we can form the
vector fields \textsc{drive }$\equiv g_{1}$, \textsc{steer
}$\equiv g_{2} $, \textsc{rotate }$\equiv $ $[g_{2},g_{1}],$ and
\textsc{slide }$\equiv \lbrack \lbrack g_{1},g_{2}],g_{1}]$,
allowing us to move anywhere in the
configuration manifold $M$. The car kinematics $\,\dot{x}%
=g_{1}u_{1}+g_{2}u_{2}\,$ is thus expanded as:
\begin{equation*}
\left(
\begin{array}{c}
\dot{x} \\
\dot{y} \\
\dot{\theta} \\
\dot{\phi}%
\end{array}
\right) =\text{\textsc{drive}}\cdot u_{1\text{
}}+\text{\textsc{steer}}\cdot u_{2\text{ }}\equiv \left(
\begin{array}{c}
\cos \theta \\
\sin \theta \\
\frac{1}{L}\tan \phi \\
0%
\end{array}
\right) \cdot u_{1\text{ }}+\left(
\begin{array}{c}
0 \\
0 \\
0 \\
1%
\end{array}
\right) \cdot u_{2\text{ }}.
\end{equation*}

The \textit{parking theorem} says: One can get out of any parking
lot that is larger than the car.

\begin{figure}[tbh]
\centerline{\includegraphics[height=4cm]{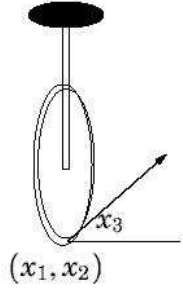}} \caption{The
unicycle.} \label{unicycle}
\end{figure}

\paragraph{The Unicycle Example.}
\label{unicyc}

Now, consider the unicycle example (see Figure \ref{unicycle}).
Here we have
\begin{eqnarray*}
g_{1} &=&\left(
\begin{array}{c}
\cos x_{3} \\
\sin x_{3} \\
0%
\end{array}%
\right) ,\qquad g_{2}=\left(
\begin{array}{c}
0 \\
0 \\
1%
\end{array}%
\right) , \\
\lbrack g_{1},g_{2}] &=&\left(
\begin{array}{c}
\sin x_{3} \\
-\cos x_{3} \\
0%
\end{array}%
\right) .
\end{eqnarray*}%
The unicycle system is full rank and therefore controllable.

\subsubsection{Controllability Condition}

Nonlinear controllability is an extension of linear
controllability. The
nonlinear SIMO system%
\begin{equation*}
\dot{x}=f(x)+g(x)\,u
\end{equation*}%
is controllable if the set of vector--fields%
\begin{equation*}
\{g,[f,g],...,[f^{n-1},g]\}
\end{equation*}%
is independent.

For example, for the kinematic kar system of the form
(\ref{Car2}), the \textit{nonlinear controllability criterion}
reads: If the Lie bracket tree:

$g_{1},$ $\ g_{2},$ $\ [g_{1},g_{2}],$ $\ [[g_{1},g_{2}],g_{1}],$
$\ [[g_{1},g_{2}],g_{2}],$ \ $[[[g_{1},g_{2}],g_{1}],g_{1}],$ \ $\
[[[g_{1},g_{2}],g_{1}],g_{2}],$ $\ [[[g_{1},g_{2}],g_{2}],g_{1}],$
$\ [[[g_{1},g_{2}],g_{2}],g_{2}],...$\newline -- has \textit{full
rank} then the system is \textit{controllable}
\cite{Isidori,Schaft,Goodwine}. In this case the combined input
\begin{equation*}
\left( u_{1},u_{2}\right) =\left\{
\begin{array}{c}
(1,0),\qquad t\in \lbrack 0,\varepsilon ] \\
(0,1),\qquad t\in \lbrack \varepsilon ,2\varepsilon ] \\
(-1,0),\qquad t\in \lbrack 2\varepsilon ,3\varepsilon ] \\
(0,-1),\qquad t\in \lbrack 3\varepsilon ,4\varepsilon ]%
\end{array}%
\right.
\end{equation*}%
gives the motion $x(4\varepsilon )=x(0)+\varepsilon ^{2}\left[ g_{1},g_{2}%
\right] +O(\varepsilon ^{3})$, with the flow given by
\begin{equation*}
F_{t}^{[g_{1},g_{2}]}=\lim_{n\rightarrow \infty }\left( F_{\sqrt{t/n}%
}^{-g_{2}}F_{\sqrt{t/n}}^{-g_{1}}F_{\sqrt{t/n}}^{g_{2}}F_{\sqrt{t/n}%
}^{g_{1}}\right) ^{n}.
\end{equation*}

\subsection{Adaptive Lie--Derivative Control}

In this subsection we develop the concept of \textit{machine
learning} in the framework of Lie--derivative control formalism
(see (\ref{liecon}) above). Consider an $n-$dimensional, SISO system in the
standard affine form (\ref{siso}), rewritten here for convenience:
\begin{equation}\label{SISO}
\dot{x}(t)=f(x)+g(x)\,u(t),   \qquad  y(t)=h(x),
\end{equation}

As already stated, the feedback control law for the system
(\ref{SISO}) can be defined using Lie derivatives $\mathcal{L}_{f}h$\ and $%
\mathcal{L}_{g}h$\ of the system's output $h$ along the
vector--fields $f$ and $g$.

If the SISO system (\ref{SISO}) is a relatively simple
(quasilinear) system
with \emph{relative degree}\footnote{%
Relative degree equals the number of differentiations of the
output function $y$ required to have the input $u$ appear
explicitly. Technically, the system (\ref{SISO}) is said to have
relative degree $r$ at the point $x^{0}$ if (see
\cite{Isidori,Schaft})
\par
(i) $\mathcal{L}_{g}\mathcal{L}_{f}^{k}h(x)=0$ for all $x$ in a
neighborhood of $x^{0}$ and all $k<r-1$, and
\par
(ii) $\mathcal{L}_{g}\mathcal{L}_{f}^{r-1}h(x^{0})\neq 0,$\newline
where $\mathcal{L}_{f}^{k}h$ denotes the $k$th Lie derivative of $h$ along $%
f$.} $=1$, it can be rewritten in a quasilinear form
\begin{equation}
\dot{x}(t)=\gamma _{i}(t)\,f_{i}(x)+d_{j}(t)\,g_{j}(x)\,u(t),
\label{reldeg1}
\end{equation}
where $\gamma _{i}$ ($i=1,...,n$)\ and $d_{j}$ ($j=1,...,m$) are
system's parameters, while $f_{i}$ and $g_{j}$ are smooth
vector--fields.

In this case the feedback control law for \emph{tracking} the
\emph{reference signal} $y_{R}=y_{R}(t)$ is defined as (see
\cite{Isidori,Schaft})
\begin{equation}
u=\frac{-\mathcal{L}_{f}h+\dot{y}_{R}+\alpha \left( y_{R}-y\right) }{%
\mathcal{L}_{g}h},  \label{ctrl}
\end{equation}%
where $\alpha $ denotes the feedback gain.

Obviously, the problem of reference signal tracking is relatively
simple and straightforward if we know all the system's parameters
$\gamma _{i}(t)$\ and $d_{j}(t)$ of (\ref{reldeg1}). The question
is can we apply a similar control law if the system parameters are
unknown?

Now we have much harder problem of \textit{adaptive signal
tracking}. However, it appears that the feedback control law can
be actually cast in a similar form (see
\cite{Sastri},\cite{Gomez}):
\begin{equation}
\widehat{u}=\frac{-\widehat{\mathcal{L}_{f}h}+\dot{y}_{R}+\alpha
\left( y_{R}-y\right) }{\widehat{\mathcal{L}_{g}h}},
\label{adaptctrl}
\end{equation}%
where Lie derivatives $\mathcal{L}_{f}h$ and $\mathcal{L}_{g}h$ of (\ref%
{ctrl}) have been replaced by their estimates
$\widehat{\mathcal{L}_{f}h}$ and $\widehat{\mathcal{L}_{g}h}$,
defined respectively as
\begin{equation*}
\widehat{\mathcal{L}_{f}h}=\widehat{\gamma _{i}}(t)\,\mathcal{L%
}_{f_{i}}h,\qquad \widehat{\mathcal{L}_{g}h}=\widehat{d_{j}}%
(t)\,\mathcal{L}_{g_{i}}h,
\end{equation*}%
in which $\widehat{\gamma _{i}}(t)$ and $\widehat{d_{j}}(t)$ are
the estimates for $\gamma _{i}(t)$\ and $d_{j}(t)$.

Therefore, we have the straightforward control law even in the
uncertain case, provided that we are able to estimate the unknown
system parameters. Probably the best known \textit{parameter
update law} is based on the so--called \textit{Lyapunov criterion}
(see \cite{Sastri}) and given by
\begin{equation}
\dot{\psi}=-\gamma \,\epsilon \,W,  \label{Lyap}
\end{equation}%
where $\psi =\{\gamma _{i}-\widehat{\gamma
_{i}},d_{j}-\widehat{d_{j}}\}$ is the parameter estimation error,
$\epsilon =y-y_{R}$ is the output error, and $\gamma $ is a
positive constant, while the matrix $W$ is defined as:
\begin{eqnarray*}
W &=&\left[ W_{1}^{T}\,W_{2}^{T}\right] ^{T},\qquad \text{with} \\
W_{1} &=&\left[
\begin{array}{c}
\mathcal{L}_{f_{1}}h \\
\vdots  \\
\mathcal{L}_{f_{n}}h%
\end{array}%
\right] ,\qquad W_{2}=\left[
\begin{array}{c}
\mathcal{L}_{g_{1}}h \\
\vdots  \\
\mathcal{L}_{g_{m}}h%
\end{array}%
\right] \cdot \frac{-\widehat{\mathcal{L}_{f}h}+\dot{y}_{R}+\alpha
\left( y_{R}-y\right) }{\widehat{\mathcal{L}_{g}h}}.
\end{eqnarray*}

The proposed adaptive control formalism
(\ref{adaptctrl}--\ref{Lyap}) can be efficiently applied wherever
we have a problem of tracking a given signal with an output of a
SISO--system (\ref{SISO}--\ref{reldeg1}) with unknown parameters.

\section{Conclusion}

In this paper we have presented two approaches to the human operator modeling: linear control theory approach and nonlinear control theory approach, based on the fixed and adaptive versions of a single-input single output Lie-Derivative controller. Our future work will focus on the generalization of the adaptive Lie-Derivative controller to MIMO systems. It would give us a rigorous closed--form model for model--free neural
networks.

\end{document}